\documentclass[11pt]{amsart}\usepackage[latin1]{inputenc}
\usepackage[T1]{fontenc}
\usepackage[english]{babel}
\usepackage{amsmath,amsthm,amssymb}
\usepackage{graphicx,epsfig,epstopdf}
\newtheorem{theo}{Theorem}[section]
\newtheorem{prop}{Proposition}[section]
\newtheorem{cor}{Corollaire}[section]
\newtheorem{lem}{Lemma}[section]
\newtheorem{defi}{Definition}[section]
\theoremstyle{remark}
\newtheorem{notat}{Notation}[section]
\newtheorem{ex}{Example}[section]
\newtheorem{rem}{Remark}

\def\A{\mathcal A}
\def\N{\mathbb{N}}
\def\Z{\mathbb{Z}}
\def\R{\mathbb{R}}
\def\T{\mathbb{T}}

\def\H{\mathcal H}
\def\FH{\mathcal{FH}}

\def\Rs{\mathcal R_{\sigma}}

\def\Rsi#1{\mathcal R_{\sigma}^{#1}}
\def\Rssi#1{\mathcal R_{S\sigma}^{#1}}

\begin{document}


\title[Product of substitutions  with the same matrix]{Random product of substitutions  with the same incidence matrix}
\author{ Pierre Arnoux}
\address{Institut de Math\'ematiques de Luminy (UPR 9016),
         163 Avenue de Luminy, case 907,
         13288 Marseille cedex 09,
         France}
\email{arnoux@iml.univ-mrs.fr}
\author{Masahiro Mizutani}
\address{Department of Management Science,
         Daito Bunka University,
         1-9-1 Takashimadaira, Itabashi-ku, Tokyo 175-8571
         JAPAN}
\email{mizutani@ic.daito.ac.jp}
\author{ Tarek Sellami}
\address{ Sfax University, Faculty of sciences of Sfax, Department of mathematics, Route Soukra
BP 802, 3018 Sfax, Tunisia.}
\email{tarek.sellami.math@gmail.com}

\subjclass[2000]{}
\date{\today}

\begin{abstract} Any infinite sequence of substitutions with the same matrix of the Pisot type defines a symbolic dynamical system which is minimal. We prove that, to any such sequence, we can associate a compact set (Rauzy fractal) by projection of the stepped line associated with an element of the symbolic system on the contracting space of the matrix. We show that this Rauzy fractal depends continuously on the sequence of substitutions, and investigate some of its properties; in some cases, this construction gives a geometric model for the symbolic dynamical system.
\end{abstract}

\maketitle

\section{Introduction}

It is well-known that any unimodular irreducible Pisot substitution defines a dynamical system that is a finite extension of a toral translation, and it is conjectured that such a dynamical system is in fact measurably equivalent to a toral translation, see \cite{PF}, \cite{GR}. A geometric model of the symbolic system can be obtained by projecting the discrete line associated to a fixed point of the substitution along its asymptotic direction, as we explain more precisely in the next section. 

In some cases, this can be extended to systems generated by an infinite sequence of substitutions belonging to a finite set $S$ (so-called $S$-adic systems), the best example being that of the sturmian sequences,  which are almost 1-1 extensions of irrational circle rotations,  the adic expansion being given by the additive continuous fraction expansion of the angle.  

It would be very interesting to be able to generalize this property; however, in the general case, we cannot expect it to hold for all sequences: already, one can check that  the closure of the set of all sturmian sequences contains periodic sequences, with bounded complexity and finite $S$-adic expansion; it is a degenerate case, where the symbolic model is finite. The existence of non balanced episturmian sequences (see \cite{JCSFLZ}) is another obstruction; in that case, the classical construction of the Rauzy fractal by projection cannot work, since the projected set is well defined, but not bounded. The existence of minimal, but not uniquely ergodic $S$-adic systems associated with interval exchange maps is a third one : for the corresponding symbolic system, the frequency is not defined, hence the discrete line (see below) associated with a symbolic sequence has no asymptotic direction, so the projection is not even defined.

In this paper, we solve the problem in a restricted case: we consider a matrix $A$ with positive integer coefficients of the Pisot type (that is, all eigenvalues are nonzero, and all eigenvalues except one are strictly smaller than 1 in modulus), and a finite set $S=\{\sigma_1, \ldots, \sigma_k\}$ of substitutions with the same matrix $A$. We first prove that any infinite sequence of elements of $S$ defines a minimal symbolic system, by defining a generalized fixed point (the {\em limit point}).  We then prove, using a generalized prefix-suffix expansion,  that this limit point stays within bounded distance of the expanding line of the matrix, and that  we can associate a compact set by projection of this limit point on the contracting space. We also prove,  by using a generalization of the classical IFS theorem of Hutchinson \cite{H},  that this set depends continuously on the sequence of substitutions, and investigate some of its properties. 

In section 2, we fix the notations, and show that a primitive sequence of substitutions on an alphabet $\A$ defines a minimal dynamical system in $\A^{\N}$. In section 3, we prove that if all the substitutions have the same  incidence matrix of the Pisot type, we can generalize to that case the projection construction of the Rauzy fractal. In section 4, we generalize the classical iterated function system to the case of an infinite sequence of contractions, and apply it to our case to prove that the generalized Rauzy fractal depends continuously on the sequence of substitutions. In section 5, we give a few examples, and make some remarks.

\section{Substitutions and adic systems}

\subsection{General setting}
Let $\A:=\{1,...,d\}$ be a finite set of cardinality $d>1$,  called the alphabet. We denote by  $\A^*$ the free monoid  on  $\A$, (set of finite words on the alphabet $\A$,  with empty word denoted by $\varepsilon$, endowed with the concatenation map). We denote by $\A^{\N}$  the set of infinite sequences on $\A$, with the natural product topology.

The length $|w|$ of a word $w \in \A^n$ with $n \in \mathbb{N}$ is defined
as $\vert w\vert = n$. For any letter $a\in\A$, we denote the number
of occurrences of $a$ in $w = w_1w_2 \ldots w_{n-1}w_n$ by $\vert w\vert_a
=\sharp \{i  | w_i = a\}$. We will denote by  $l : \A^* \mapsto \mathbb{N}^d : w\mapsto(\vert w\vert
_a)_{a\in\A}\in\mathbb{N}^d$  the natural homomorphism (abelianization map) obtained by
abelianization of the free monoid.

A substitution over the alphabet $\A$ is a nonerasing endomorphism of the free
monoid $\A^*$. To any substitution $\sigma$, one can associate its
incidence matrix $M$,  which is the $d\times d $ matrix obtained by
abelianization, i.e. $M_{i,j} = \vert \sigma(j) \vert_i$. By construction, one has
$l(\sigma(w))=Ml(w)$ for any word $w\in\A^*$.\

\begin{defi} A substitution $\sigma$ is {\em primitive } if there exists an integer
$k$ such that, for each pair $(a, b)\in \A^2$, $\vert \sigma^k(a)
\vert_b > 0$. 
\end{defi}

It is equivalent to suppose that the incidence matrix is primitive, that is, this matrix has a strictly positive power. We will always suppose that the substitution is primitive, this
implies that for all letter $j \in \A$ the length of the successive iterations
$\sigma^k(j)$ tends to infinity.

\begin{rem} Since the incidence matrix of a primitive substitution is a
primitive matrix, by the Perron-Frobenius theorem, it has a simple real
positive dominant eigenvalue $\beta$.
\end{rem}

\subsection{Dynamical system defined by a primitive substitution}

A substitution $\sigma$ naturally extends to the set $\A^{\N}$ of infinite sequences with value in $\A$, by defining $\sigma(u_1 u_2\ldots)=\sigma(u_1)\sigma(u_2)\ldots$. We say that a sequence $u\in\A^{\N}$ is a periodic point of $\sigma$ if there exists some integer $k$ such that $\sigma^k(u)=u$. One easily proves that any primitive substitution has periodic points, and that two periodic points with the same initial letters are equal, hence there are at most $d$ periodic points.

We denote by $S$ the shift on $\A^{\N}$, defined by $S(u)=v$, where $v$ is the sequence such that, for all $i\in\N$, $v_i=u_{i+1}$.

\begin{defi}
Let $\sigma$ be a primitive substitution, and let $u$ be a periodic point of $\sigma$. Let $\Omega_{\sigma}=\overline{\{S^n(u)|n\in\N\}}$ be the closure of the orbit of $u$ by the shift.

The dynamical system defined by the substitution $\sigma$ is $(\Omega_{\sigma},S)$.
\end{defi}

\begin{rem} The primitivity condition implies that any word occuring in a periodic word of $\sigma$ also appears in any other periodic word of $\sigma$, hence all these periodic words have same closure for the shift, and the set $\Omega_{\sigma}$ does not depend on the particular periodic word we have chosen, only on the substitution $\sigma$.
\end{rem}

\subsection{Pisot substitutions} 

We recall a classical definition of algebraic number theory:  
\begin{defi}
A Pisot number is an algebraic integer $\beta>1$ such that each Galois
conjugate $\beta^{(i)}$ of $\beta$ satisfies $\mid\beta^{(i)}\mid<1$.
\end{defi}

By analogy, we define substitutions of Pisot type and irreducible substitutions of Pisot type: 

\begin{defi} A substitution is of {\em Pisot type} if the dominant eigenvalue of its incidence matrix is a Pisot number. A substitution of Pisot type  is {\em (algebraically) irreducible} if the characteristic polynomial of its incidence matrix is irreducible over the field of the rational numbers.
\end{defi}

\begin{rem} It is equivalent to say that all the eigenvalues of the incidence matrix are not equal to 0, and all except one are strictly smaller than 1 in modulus. Any irreducible substitution of Pisot type is primitive (see  \cite{VCAS}), but the converse is not true: there are primitive substitutions of Pisot type which are reducible.
\end{rem}

\begin{rem}
The fixed point of a substitution of Pisot type can never be periodic. 
\end{rem}
\subsection{$S$-adic systems}

Dynamical system defined by substitutions form a small class; in particular, there is only a countable number of such system. For example, the only circle rotations whose natural coding is given by a subsitution are those whose rotation number is a quadratic integer of a special type (Galois integer).

It is however possible to enlarge this class by replacing one iterated substitution by an infinite sequence of substitutions chosen in a finite set $S$; this can be considered as a kind of continued fraction expansion, the substitution dynamical systems corresponding to  periodic continued fraction expansions. We obtain in this way the class of $S$-adic systems, but some care must be taken to generalize the primitivity condition.

\begin{defi}Let $S$ be a finite set of substitutions on a fixed alphabet $\A$. Let $(\sigma_n)_{n\in\N}$ be a sequence of substitutions in $S$. We say that the sequence is {\em primitive} if, for any $n\in\N$, there is an integer $k$ such that the substitution $\sigma_n\sigma_{n+1}\ldots\sigma_{n+k}$ has a strictly positive incidence matrix. 
\end{defi}

We can extend the notion of periodic points to such a sequence of substitutions.

\begin{defi} Let $(\sigma_n)_{n\in\N}$ be a sequence of substitutions in $S$. We say that  $u\in\A^{\N}$ is a {\em limit point} of the sequence $(\sigma_n)_{n\in\N}$ if there exists a sequence $(u^{(n)})$ of points of $\A^{\N}$ such that $u=u^{(0)}$ and $u^{(n)}=\sigma_n(u^{(n+1)})$.
\end{defi}

\begin{prop} Any primitive sequence of substitutions has a finite and nonzero number of limit points.

Any finite word that occurs in a limit point of a primitive sequence of substitutions occurs in all the other limit points.
\end{prop}
\begin{proof}
With  any substitution $\sigma$, we associate the map $f: \A\to \A$ which sends any letter $a\in\A$ to the first letter of  $\sigma(a)$. Given a sequence $(\sigma_n)$ of substitutions, we consider the sequence
$f_1(f_2\ldots (f_n(\A))\ldots))$ of subsets of $\A$. This is a decreasing sequence of nonempty subsets of $\A$, so their intersection is nonempty. Let $a$ be any element of this intersection; by construction, there exists a sequence $(a_n)$ of elements of $\A$ such that $a_0=a$ and $a_n=f_n(a_{n+1})$ for all $n$.

Consider the sequence of words $U_n=\sigma_0\sigma_1\ldots\sigma_{n-1}(a_{n})$. Since $a_n$ is the first letter of $\sigma_n(a_{n+1})$, $U_n$ is by construction a prefix of $U_{n+1}$. The primitivity condition implies that the length of $U_n$ tends to infinity, since for any $n$ there is a $k$ such that   $\sigma_n\sigma_{n+1}\ldots\sigma_{n+k}$ sends any letter to a word of length at least $d>1$.

This sequence of finite words defines a unique infinite sequence $u$ which admits these words for prefix. By shifting the sequences $(\sigma_n)$ and $(a_n)$, we define in the same way points $u^{(n)}$ which, by construction, satisfy $u^{(n)}=\sigma_n(u^{(n+1)})$. Hence any primitive sequence of substitutions has a limit point.  

Remark that any limit point arises in this way: if we have a sequence $(u^{(n)})$ of points of $\A^{\N}$ such that $u^{(n)}=\sigma_n(u^{(n+1)})$, and if we define $a_n$ to be the first letter of $u^{(n)}$, the words $\sigma_0\sigma_1\ldots\sigma_{n-1}(a_{n})$ are by definition prefixes of $u^{(0)}$, whose lengths tend to infinity.

Since, for any sequence $(a_n)$ such that  $a_n=f_n(a_{n+1})$, $a_n$ defines $a_p$ for all $p<n$, there can be at most $d$ such sequences (because there are at most $d$ prefixes of length $n$ of such sequences for any $n$). Hence there are at most $d$ limit points.

Let $u$ and $v$ be two limit points, and let $(a_n)$ and $(b_n)$ be the corresponding sequences of letters. Let $U$ be a finite word which occurs in $u$. Then there exists an $n$ such that $U$ occurs in the word $\sigma_0\sigma_1\ldots\sigma_{n-1}(a_{n})$. By primitivity, we can find $k$ such that the letter $a_n$ occurs in the word $\sigma_n\ldots\sigma_{n+k}(b_{n+k})$. Hence $U$ occurs in the word $\sigma_0\sigma_1\ldots\sigma_{n+k}(b_{n+k})$, which is a prefix of $v$. This proves that a word which occurs in a limit word occurs in all the other limit words.
\end{proof}
\begin{rem}
The special case of a constant sequence corresponds to the fixed point of a substitution; in that case, the definition of primitivity of a sequence reduces to the definition of primitivity of a substitution.
\end{rem}

\begin{defi} Let $S$ be a finite set of substitutions, and let $\sigma=(\sigma_n)_{n\in\N}$ be a primitive sequence of substitutions. The $S$-adic system generated by the sequence of substitutions $\sigma$ is the dynamical system $(\Omega_{\sigma},S)$, where $\Omega_{\sigma}$ is the closure of the orbit by the shift of a limit point of $\sigma$.
\end{defi}

The previous proposition shows that the set $\Omega_{\sigma}$ does not depend on the choice of the limit point $u$, but only on the sequence $\sigma$.

\begin{prop}The $S$-adic system generated by a primitive sequence of substitutions is minimal.
\end{prop}

\begin{proof} It is enough to prove that any word which occurs in a point of the system occurs with bounded gap, and it is enough to prove this for a limit point. Using the same notations as above, let $U$ be the finite word under consideration. Then we can find $n$ such that $U$ occurs in   $\sigma_0\sigma_1\ldots\sigma_{n-1}(a_{n})$, and by primitivity we can find $k$ such that $U$ occurs in     
$\sigma_0\ldots\sigma_{n+k}(b)$ for any $b$ in $\A$. Let $L$ be the maximum length of the words 
$\sigma_0\ldots\sigma_{n+k}(b)$; since the limit point can be decomposed in these words, the word $U$ occurs in the limit point with gaps bounded by $2L$.
\end{proof}

\begin{ex} A very simple example of a family of $S$-adic system is given by the sturmian sequences, or rotation sequences, which are almost one-to-one extension of rotations. Here, the alphabet $\A=\{0,1\}$, and the set $S$ of substitutions is $S=\{\sigma_0, \sigma_1\}$, where $\sigma_0$ is defined by $\sigma_0(0)=0 , \quad \sigma_0(1)=10$ and $\sigma_1$ is defined by $\sigma_1(0)=01 , \quad \sigma_1(1)=1$. In that case, a sequence is primitive if and only if it is not eventually constant, and the sequence of substitutions which defines a given sturmian system is given by  the additive  continued fraction expansion of the angle of the associated rotation.
\end{ex}

\begin{ex} A more tricky example is given by the natural symbolic dynamics of interval exchange maps. Here also we get a finite set of substitutions, associated with the Rauzy induction which generalizes the classical continued fraction. This example shows that we cannot expect to generalize to the $S$-adic case all the properties of substitution dynamical systems; in particular, there are $S$-adic systems defined by a primitive sequence of substitutions which are not uniquely ergodic, as shown by the existence of minimal but not uniquely ergodic interval exchange maps. 
\end{ex}

\section{Geometric models for substitution dynamical systems and $S$-adic systems}

Substitution dynamical systems are one of the simplest families of dynamical systems, and they are present in every framework where self-similarity occurs. It is natural to study their dynamical properties; as we have seen above, they are minimal, and a simple application of the Perron-Frobenius theorem shows that any letter, and any finite word, has a well defined frequency, which proves that the system is uniquely ergodic.

Hence, a substitution dynamical system is a measured dynamical system in a natural way, and we can study its ergodic properties, look for eigenvalues, and look for the maximal equicontinuous factor. It is well-known that the Pisot condition plays a role here: the work of Rauzy and other people (\cite{GR}, \cite {PF} has shown that any unimodular irreducible substitution of Pisot type on $d$ letters is a finite extension of a translation of the torus $\T^{d-1}$, and it is conjectured that it is an almost one-to-one extension. 

We will briefly review  a geometric construction for this maximal equicontinuous factor of a unimodular Pisot substitution $\sigma$. We first need to introduce some algebraic formalism in order to embed the periodic point $u$ in a linear subspace spanned by the algebric conjugates of the dominant eigenvalue of the incidence matrix of $\sigma$; the closure of the "projections" of the prefixes of $u$ will form the so-called central tile or Rauzy fractal.

\subsection{Rauzy fractal associated to a unimodular Pisot substitution}
\begin{defi}
A stepped line $L=(x_n)$ in $\mathbb{R}^d$ is a sequence (finite or infinite) of points in $\mathbb{R}^d$
such that $x_{n+1}-x_n$ belongs to a finite set.

A canonical stepped line is a stepped line such that $x_0=0$ and for all $n\geq 0$, $x_{n+1}-x_n$ belongs to the canonical basis of $\mathbb{R}^d$.
\end{defi}
Using the abelianization map, to any finite or infinite word $W$, we can
associate a canonical stepped line in  $\mathbb{R}^d$ as the sequence $(l(P_n))$,
where $P_n$ is the prefix of length $n$ of $W$.

The Perron-Frobenius Theorem implies that the canonical stepped line associated with a periodic point of a  primitive substitution has an asymptotic direction. If this substitution is irreducible and of Pisot type, this canonical stepped line  remains within bounded distance from the expanding line (given by the Perron-Frobenius  right eigenvector of $M_{\sigma}$). These two properties (bounded distance and Pisot type) are equivalent for algebraically irreducible substitutions.

\begin{notat}Let  $\sigma$ be an irreducible Pisot substitution. We denote by $E_s$ the stable space (or contracting space) and $E_u$ the unstable space (or expanding line) of $M_{\sigma}$.
We denote by $\pi_s$ the linear projection on $E_s$, parallel to $E_u$. 
\end{notat}
Since the canonical stepped line stays within bounded distance of $E_u$, its projection is a  bounded set in the $(d-1)$-dimensional vector space $E_s$.

\begin{defi}
Let $\sigma$ an be irreducible Pisot substitution. The Rauzy fractal (or central tile) associated with
$\sigma$ is the closure of the projection of the canonical stepped line
associated to any periodic point of $\sigma$ in the contracting plane parallel to
the expanding direction. We will denote it by $\Rs$.
\end{defi}

The Rauzy fractal can be decomposed in subtiles depending on the letter associated to
the vertex of the stepped line that is projected. For any $i\in\A$, we define  : $ \Rsi{i} := \overline{\{\pi_s(l(u_0\ldots u_{k-1}), k\in{\mathbb{N}}, u_k=i\}} $. We have obviously $\Rs=\bigcup_{i\in\A}\Rsi{i}.$

\textbf{Remark:} One can prove that the definition of $\Rs$ and $\Rsi{i}$  for $i\in\A$ does not depend on the choice of the periodic point $u\in\A$.

We will briefly recall some properties of the Rauzy fractal. If we denote by $\Gamma$ the projection by $\pi_s$ of the diagonal lattice of $\Z^d$ (that is, the set of $(n_1, \ldots, n_d)$ such that $\sum n_i=0$), the translate by $\Gamma$ of the Rauzy fractal cover the stable plane. Hence the Rauzy fractal has positive measure. The projection from the orbit of the periodic point to the Rauzy fractal extends by continuity to all of $\Omega_{\sigma}$; if we quotient by $\Gamma$, this projection gives a semi-conjugacy between $(\Omega_{\sigma},S)$ and the translation by $\pi_s(e)$ on $\Rs/\Gamma$, where $e$ is any vector in the canonical base (all $\pi_s(e)$ are equivalent modulo $\Gamma$).

It is conjectured that the semi-conjugacy is a conjugacy, $\Rs$ is a fundamental domain for $\Gamma$ and the $\Rsi{i}$ are a partition of $\Rs$, up to a set of measure 0, and this can be proved in many particular cases.

Finally, the subsets of the Rauzy fractal are solution of a set equation: 

\begin{prop}
Let $\sigma$ be a unimodular irreducible substitution of Pisot type, with incidence matrix $M$.  Let us denote by $l_i=|\sigma(i)|$ the length of $\sigma(i)$, and $\sigma(i)=w^{(i)}_1\ldots w^{(i)}_{l_i}=P^{(i)}_kw^{(i)}_k S^{(i)}_k$, where $P^{(i)}_k$ is the prefix of length $k-1$ of $\sigma(i)$ and $S^{(i)}_k$ is the suffix of length $l_i-k$ of $\sigma(i)$. 

Then the subsets $\Rsi{i}$ of the Rauzy fractal satisfy the following set equation:
$$\Rsi{i}=\bigcup_{j,k; w^{j}_k=i} M. \Rsi{j}+\pi_s(l(P^{(j)}_k))$$
\end{prop}

\subsection{Rauzy fractal associated to a sequence of substitutions}

We want to generalize this construction to the $S$-adic case. To do this, we can replace the periodic point by a limit point, but we need to have an asymptotic direction, and a bounded distance property. The most general condition under which this is possible is very unclear, and we will restrict to a very special  case.

We consider a positive integral primitive matrix $M$ with determinant $\pm1$ satisfying the Pisot condition: all the eigenvalues different from the Perron-Frobenius eigenvalue are strictly smaller than 1 in modulus. Let $S$ be a set of substitutions, all of them having $M$ as incidence matrix.

It is then clear by construction that any sequence in $S^{\N}$ is a primitive sequence of substitutions. To such a sequence, one can associate a limit point, as we have shown before. This limit point has an asymptotic direction, which is the unstable direction of the matrix $M$. 

We need to fix some notation. $\sigma=(\sigma_n)_{n\in\N}\in S^{\N}$ is a sequence of substitutions; as above, we write,  for any substitution $\sigma_0\in S$, $l_i=|\sigma_0(i)|$ (it only depends on $M$, not $\sigma_0$), and $\sigma_0(i)=w^{\sigma_0,i}_1\ldots w^{\sigma_0,i}_{l_i}=P^{\sigma_0,i}_kw^{\sigma_0,i}_k S^{\sigma_0,i}_k$,  where $P^{\sigma_0,i}_k$ is the prefix of length $k-1$ of $\sigma_0(i)$.

\begin{lem}
Let $u$ be the limit point of a sequence $\sigma\in S^{\N}$. The canonical stepped line of $u$ stays within bounded distance of the expanding line of $M$.
\end{lem}

\begin{proof} Let $U$ be a prefix of the limit point. By construction, we can write $U=\sigma_0(U^{(1)}) P_0$, where $U^{(1)}$ is some prefix of a limit point of the shifted sequence $S\sigma$, and $P_0$ is a prefix of $\sigma_0(i)$ for some $i\in\A$. We can iterate on $U^{(1)}$, and a recurrence proves that there is some $k$ such that we can write 
$$U=\sigma_0\sigma_1\ldots\sigma_{k-1}(P_k)\sigma_0\sigma_1\ldots\sigma_{k-2}(P_{k-1})\ldots\sigma_0(P_1)P_0$$
where, for all $j\le k$,  $P_j$ is a prefix of $\sigma_j(i)$, for some letter $i\in\A$.

Taking abelianization, since all substitutions have the same incidence matrix,  we get

\(\begin{array}{*{3}{c}}
$$l(U)$$ & 
 $$ = $$ & $$l( \sigma_0\sigma_1\ldots\sigma_{k-1}(P_k)\sigma_0\sigma_1\ldots\sigma_{k-2}(P_{k-1})\ldots\sigma_0(P_1)P_0)$$\\
$${}$$ & $$ = $$ & $$ l( \sigma_0\sigma_1\ldots\sigma_{k-1}(P_k)) + l(\sigma_0\sigma_1\ldots\sigma_{k-2}(P_{k-1})) +\ldots + l(P_0)$$\\
$${}$$ & $$ = $$ & $$M^kl(P_k) + M^{k-1}l(P_{k-1}+\ldots Ml(P_1) + Ml(P_1)$$\\
$${}$$& $$ = $$ & $$\sum _{j=0}^k M^j l(P_j)$$\\

\end{array}\)

 Since the stable and unstable space are invariant under $M$, the projection $\pi_s$ commutes with $M$ and we get:
$$\pi_s l(U)=\sum _{j=0}^k M^j \pi_s l(P_j)$$
One can always define an adapted  norm on the contracting space $E_s$ such that  $|M X |  \le \lambda  | X| $, for some $\lambda<1$. We obtain :

\(\begin{array}{*{3}{c}}
$$ |\pi_s(l(U))|$$ & $$ = $$ & $$|   \sum _{j=0}^k M^j \pi_s l(P_j) |$$\\
$${}$$ & $$\leq$$ & $$\sum _{j=0}^k | M^j \pi_s l(P_j) | $$\\
$$ {} $$ & $$\leq $$& $$\sum_{j=0}^k \lambda^j |\pi_s l(P_j) |  $$\\

\end{array}\)

But there is a finite number of substitutions, hence a finite number of prefixes, so that $|\pi_s l(P_j)|$ is bounded by some constant $C$; hence the finite sum for $|\pi_s l(P_j)|$ is bounded, independently  of the length of $U$, by the geometric series $\sum _{j=0}^{\infty} C\lambda^j=\frac C{1-\lambda}$.

This proves that the canonical stepped line associated with the limit point stays within bounded distance of the expanding line, and that its projection by $\pi_s$ is bounded.
\end{proof}

\begin{defi} Let $M$ be a positive integral unimodular matrix of the Pisot type. Let $\sigma$ be a sequence of substitutions with incidence matrix $M$. The Rauzy fractal associated to the sequence $\sigma$ is the closure of the projection to the stable space of the canonical line associated with a  limit point of $\sigma$. It is denoted by $\Rs$.
\end{defi}

We can as above define subsets $\Rsi i$ of $\Rs$ corresponding to a particular letter; they satisfy a set equation, which relates the Rauzy fractal associated with $\sigma$ to the Rauzy fractal associated with $S\sigma$.
\begin{prop}
We have: 
$$\Rsi{i}=\bigcup_{j,k; w^{\sigma_0,j}_k=i} M. \Rssi{j}+\pi_s(l(P^{(\sigma_0,j)}_k))$$
\end{prop}
\begin{proof}
This is an immediate consequence of the equality given above for prefixes of the limit point, $U=\sigma_0(U^{(1)}) P_0$, by projecting to the stable space and taking closure.
\end{proof}
     
We will use this set equation in the next section to prove that $\Rs$ depends continuously on the sequence $\sigma$.

\section{Iterated function systems}

\subsection{Hutchinson's theorem}
Many fractals are made up of parts that are similar to the whole. For example, the middle third Cantor set is the union of two similar copies of itself, and the von Koch curve is made up of four similar copies. These self-similarities are not only properties of the fractals: they may actually be used to define them. Iterated function system do this in a unified way and, moreover, often lead to a simple way of finding dimensions. 

A mapping $S  : \R^n\rightarrow \R^n$ is called a contraction if there is a number $c$ with $0<c<1$ such that  $| S(x)-S(y)|\le c|x-y|$ for all $x,y\in D$. Clearly any contraction is continuous. If equality holds, i.e. if  $| S(x)-S(y)| = c|x-y|$, then $S$ transforms sets into geometrically similar sets, and we call $S$ a contracting similarity.
\begin{defi}
A finite family of contractions $\{f_1,f_2,\ldots, f_m\}$, with $m\le 2$, is called an iterated function system or IFS. We call a non-empty compact subset $K$ of $\R^n$ an attractor for the IFS if $$ K = \bigcup_{i=1}^{m} f_i(K).$$
\end{defi}

The fundamental property of an iterated function system is that it determines a unique attractor, which is usually a fractal. The proof is a simple consequence of Banach fixed point theorem.

We denote by $\mathcal{H}$ the set of nonempty compact subsets  of $\R^n$, and consider the map $\Phi$ defined by: 
\begin{eqnarray*}
\Phi : &\mathcal{H}&\to\mathcal{H}\\
&K &\mapsto \bigcup_{i=1}^m f_i(K)
\end{eqnarray*}

Let $K$ a non empty  compact set from $\mathbb{R}^n$, and $x\in\mathbb{R}^n$.  We define the distance $d(x,K) = \inf_{y\in K} d(x,y)$.
We define the Hausdorff  distance between elements of $\mathcal{H}$ by $$d_\mathcal{H} (K, K') = \max (\sup_{x\in K}d(x,K), \sup_{x\in K'} d(x,K')).$$

\begin{lem}
$(\mathcal{H}, d_{\mathcal{H}})$ is a complete space.
\end{lem}

\begin{proof}
See \cite{KF}
\end{proof}

\begin{lem}
 $\Phi$ is a contracting map.
\end{lem}

\begin{proof} If $f$ is a contraction of constant $c$, then $d(f(K), f(K'))\le c d(K,K')$. If $(A_1, \ldots A_m)$ and $(B_1, \ldots, B_m)$ are finite collections of elements of $\H$, an immediate computation shows that 
$d(\cup_i A_i, \cup_i B_i)\le \max_i d(A_i,B_i)$. Hence, if all the $f_i$ are contractions on $\R^n$ of ratio less than $c$, $\Phi$ is a contraction on $\H$ of ratio less than $c$.
\end{proof}

\begin{theo}
$\Phi$ admit a  unique  fixed point.
\end{theo}

\begin{proof}
$\Phi$ is a contraction of ration $c<1$ on a complete space, hence it has a unique fixed  point by the Banach fixed point theorem.
\end{proof}

\subsection{Graph-directed Iterated Function Systems}

An IFS acts on a unique set $K$, by replacing it   by its image $ \bigcup_{i=1}^{m} f_i(K)$, and its fixed point is defined by a unique equation. In many cases,  the fixed point has a finite number of components, and is define by a set of equations, as we saw for the Rauzy fractal. This more general framework gives rise to what is called  {\em Graph-directed Iterated Function Systems} (GIFS for short).

\begin{defi}
Let G be a finite directed graph with set $V$ of vertices $\{1,\ldots,d\}$ and set
of edges $E$. Denote the set of edges leading from $i$ to $j$ by $E_{ij}$. To
each $e\in E$ associated a contractive mapping $f_e:\R^n\rightarrow
\R^n$. If for each $i$ there is some outgoing edge we call
$(G,\{f_e\})$ a GIFS.
\end{defi}

Such a GIFS defines a map $\Phi$ on $\H^d$ by 

\begin{eqnarray*}
\Phi : &\H^d&\to\H^d\\
&(K_i)_{i=1..d} &\mapsto \left(\bigcup_{j=1}^d\bigcup_{e\in E_{i,j}} f_e(K_j)  \right)_{i=1..d}
\end{eqnarray*}

If we define a distance on $\H^d$ by $d((K_i), (K'_i))=\max_i( d(K_i, K'_i))$, $\H^d$ is complete for this distance, and $\Phi$ is a contraction. The same fixed point argument shows that to a GIFS $(G,\{f_e\}_{e\in E})$ there corresponds a unique collection of nonempty compacts sets $K_1,\ldots, K_q \subset \mathbb{R}^n$ having the property that $$K_i=\bigcup_{j=1}^q \bigcup_{e\in E_{ij}} f_e(K_j).$$The collection $K_1,\ldots, K_q$ is called GIFS attractor or solution of the GIFS. 

\subsection{Fixed point associated to a sequence of IFS}

We consider now a finite family of  $k$ GIFS acting on the same set of vertices $V=\{1,\ldots, d\}$.  To each GIFS we associate the contraction $\Phi_j$ of ratio $c_j<1$ acting on $\H^d$, as defined above, with $j\in\{1, \ldots, k\}$.

We define $C=\{1,\ldots, k\}^{\N}$, the set of sequences with value in $\{1, \ldots, k\}$, with the shift acting on $C$ in a natural way. Our goal is to associate to each sequence $\epsilon=\epsilon_0, \epsilon_1,\ldots \in C$ a family of $d$ compact sets $K_{\epsilon}=(K_{\epsilon,1}, \ldots, K_{\epsilon,d})$ such that $K_{\epsilon}=\Phi_{\epsilon_0}K_{S\epsilon}$.

We consider the set  $\FH$ of continuous functions from $C$ to $\H^d$ with the uniform distance, defined for any elements  $F_1, F_2\in\FH$ by 

$$d(F_1,F_2) = \sup_{\varepsilon\in C} d_\mathcal{H}(F_1(\varepsilon), F_2(\varepsilon)).$$

\begin{lem}
$\FH$ is a complete space.
\end{lem}
\begin{proof}
Let $F_n$ a Cauchy sequences from $\FH.$ This mean for all $\alpha >0$ , we can find $N\in\mathbb{N} $ such that, for all
$ p,q>N$,  $d(F_p, F_q)< \alpha$.  Hence in particular, for any $\varepsilon\in C$,  $F_n(\varepsilon)$ is a Cauchy sequences in $\H$. 
Since $\H^d$ is a complete space,  it converge to some limit $F_{\infty}(\varepsilon)$. Hence  the sequence $F_n$ converges pointwise to some function $F_{\infty}$. 

Taking $\alpha$ and $N$ as above, we see that, for any $n>N$, $d_\mathcal{H}(F_n(\varepsilon), F_{\infty}(\varepsilon))<\alpha$, hence $d(F_n,F_{\infty})<\alpha$. This proves that the convergence of $F_n$ to $F_{\infty}$ is in fact uniform, hence $F_{\infty}$ is continuous, and the space $\FH$ is complete.
\end{proof}

\begin{lem}
We consider the map $\Phi : \FH \longrightarrow \FH$ defined by :  $\Phi(F) = G$ such that $G(\varepsilon) = \Phi_{\varepsilon_0}(F(S\varepsilon))$. Then $\Phi$ is a contracting map of ratio $c<1$.
\end{lem}

\begin{proof}
Let $F_1$ and $F_2$ in $\mathcal{FH}$. 

We know that, for a given $\epsilon$,  $\Phi_{\epsilon_0}$ is a contracting map, hence there exist $c_{\epsilon_0}<1$, depending only on  ${\epsilon_0}$, such that $$d_\mathcal{H} (\Phi_{\epsilon_0}F_1(S\epsilon), \Phi_{\epsilon_0}F_2(S\epsilon))\leq c d_\mathcal{H} (F_1(S\epsilon), F_2(S\epsilon))$$

If we define $c=\sup_{i=1..k} c _{i}$, we obtain 

$$ \sup_{\epsilon\in C}d_\mathcal{H}(\Phi F_1(\epsilon),\Phi F_2(\epsilon)) \leq c \sup_{\epsilon\in C}  d_\mathcal{H} (F_1(S\epsilon), F_2(S\epsilon))    $$

Hence $d(\Phi(F_1), \Phi(F_2))\leq c d(F_1,F_2)$; since $c<1$, we have proved that $\Phi$ is a contraction.
\end{proof}

\begin{theo}
$\Phi$ admet a  unique fixed point.
\end{theo}

\begin{proof}
Same as above: $\Phi$ is a strictly contracting map in a complete space $\FH$, hence it  admits a unique fixed point.
\end{proof}

\subsection{Application to generalized Rauzy fractals}

We proved in the previous section that, if we have a finite set $S$ of substitution with the same matrix of the Pisot type, we can associate to any sequence $\sigma$ of substitutions in $S$ a generalized Rauzy fractal $\Rs$, and that this Rauzy fractal satisfies the following equation: 

$$\Rsi{i}=\bigcup_{j,k; w^{\sigma_0,j}_k=i} M. \Rssi{j}+\pi_s(l(P^{(\sigma_0,j)}_k))$$

We can view $(\Rsi{i})$ as a map $\mathcal R: S^{\N}\to \H^d$, $\sigma\mapsto (\Rsi{i})$, where $\H$ is the set of nonempty compact subsets of the contracting plane of $M$. Each substitution $\sigma_0\in S$ defines a contraction $\Phi_{\sigma_0}$ on $\H^d$, the contraction ratio being given by the contraction ration of the restriction of $M$ to the contracting space.

We are exactly in the conditions to apply the theorem above, and we obtain: 

\begin{cor}
Let $S$ be a finite set of substitutions with the same unimodular Pisot incidence matrix $M$. To any sequence $\sigma\in S^{\N}$, we can associate a Rauzy fractal $\Rs=\cup_{i=1}^d\Rsi{i}$ in the contracting space of $M$. This Rauzy fractal  satisfies the relation 
$$\Rsi{i}=\bigcup_{j,k; w^{\sigma_0,j}_k=i} M. \Rssi{j}+\pi_s(l(P^{(\sigma_0,j)}_k))$$
It depends continuously on $\sigma$ for the natural product topology on $S^{\N}$.
\end{cor}
\section{Exemples}

A classic example is given by  the Tribonacci
substitution and the flipped Tribonacci substitution, i.e.,

\begin{center}
$\sigma_1:\left\{
\begin{array}{ll}
a\rightarrow ab\\
b\rightarrow ac\\
c\rightarrow a
\end{array}\right.$
\hspace{1cm}and \hspace{1cm} $\sigma_2:\left\{
\begin{array}{ll}
a\rightarrow ab\\
b\rightarrow ca\\
c\rightarrow a
\end{array}\right.$
\end{center}
\vspace{1cm} The incidence matrix of $\sigma_1$ and $\sigma_2$ is
$\begin{pmatrix}
1 & 1 & 1 \\
1 & 0 &0 \\
0 & 1 & 0 \\
\end{pmatrix}
$. The dominant eigenvalue satisfies the relation $X^3-X^2-X-1=0$, hence the
name Tribonacci for the substitution.

 The Rauzy fractal of the first substitution is a topological disc \cite{PA},
simply connected , while it is a well known fact that the second fractal is not
simply connected, and has infinitely generated fundamental group, see \cite{ASJT}; compare Figure[1]. \vspace{0.5cm}
\begin{figure}[h]
\begin{center}
\label{compare_fractals}
\scalebox{0.4}{\includegraphics{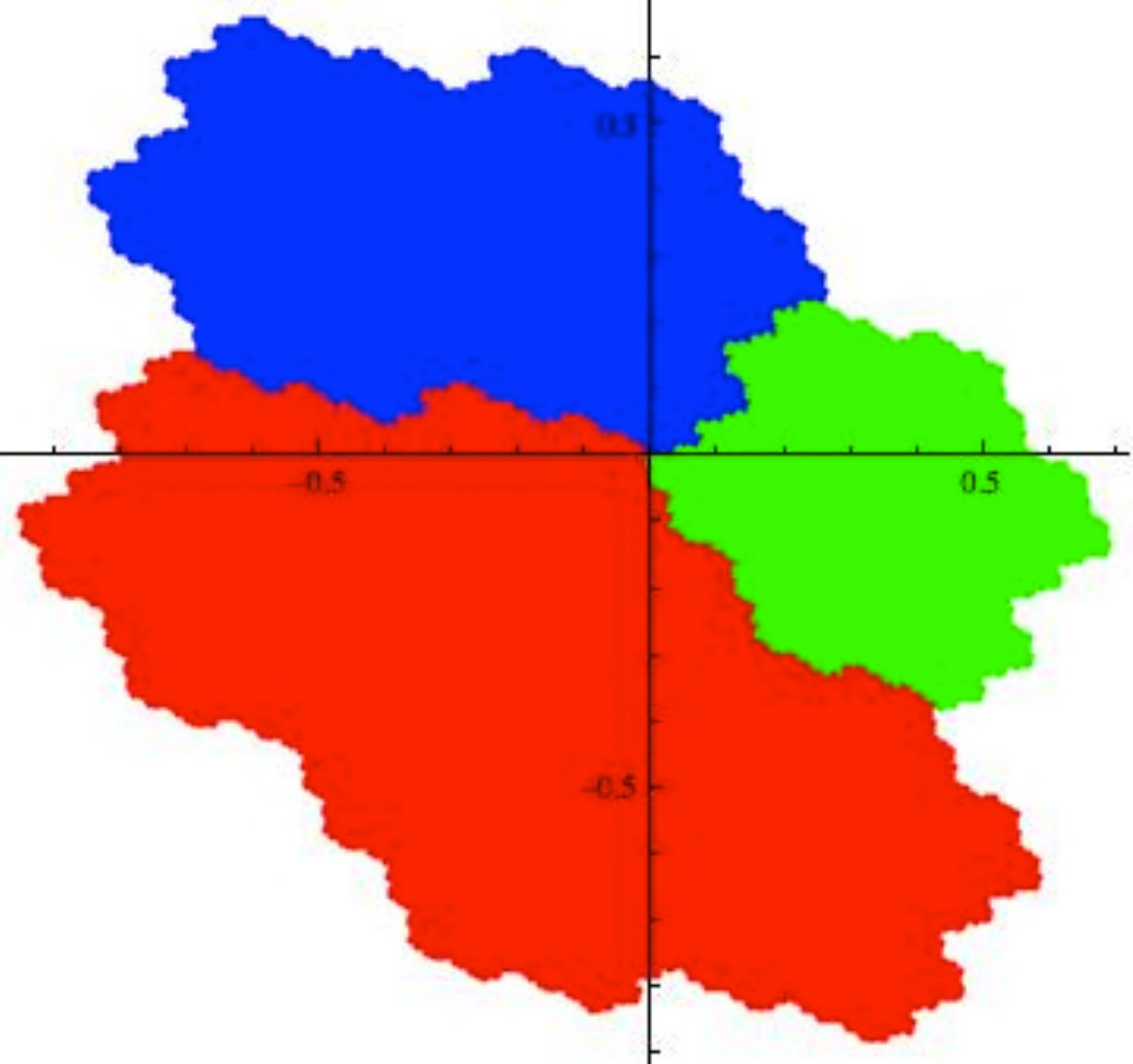}}
\vspace{0.5cm}
\scalebox{0.4}{\includegraphics{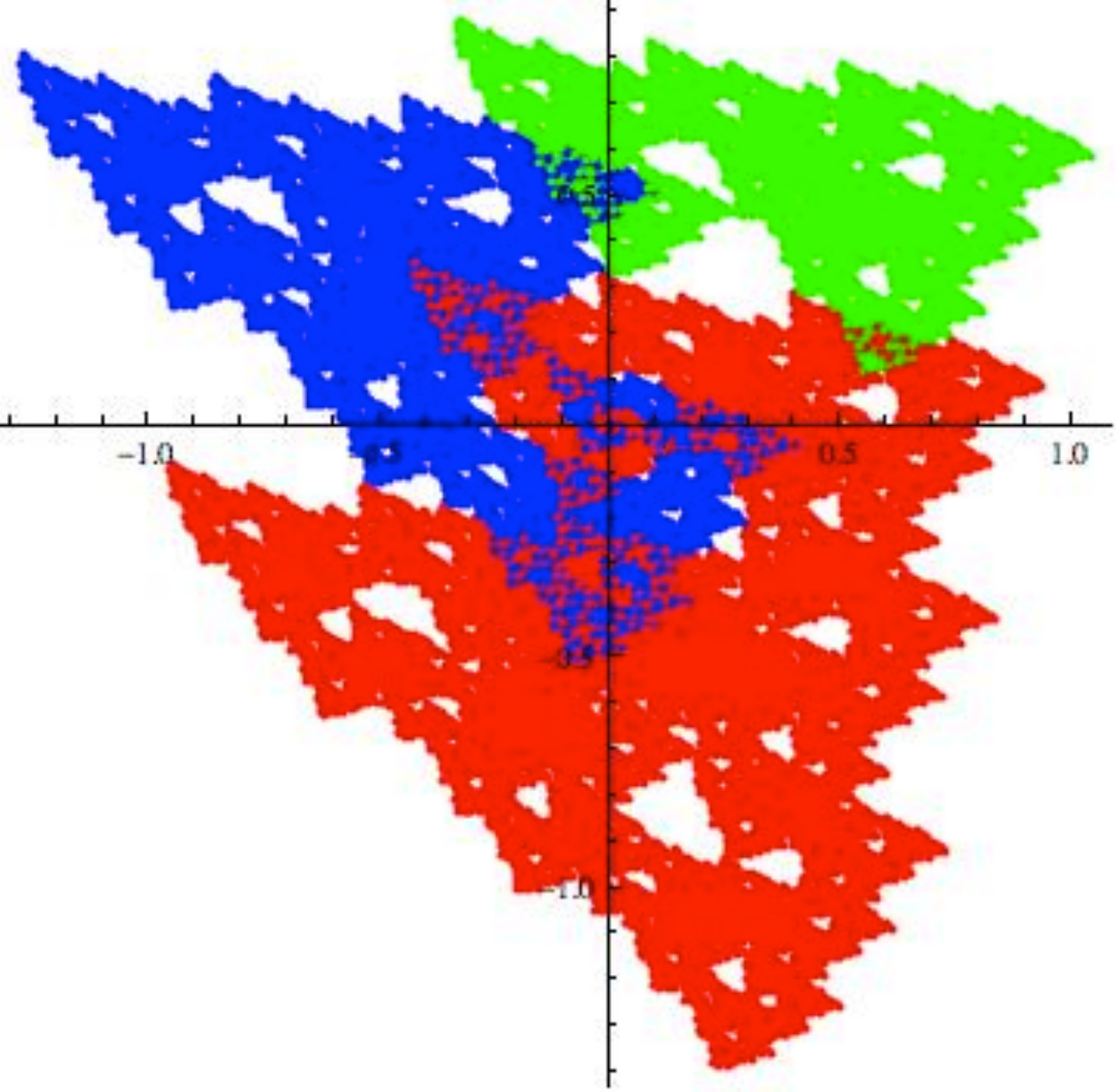}}
\caption{ The Rauzy  fractals of  $\sigma_1$ and $\sigma_2$}
\end{center}
\end{figure}
We obtain some nice Rauzy fractals for different sequences $w_n$, in particular we can find some examples with a nontrivial but finitely generated fundamental group.

\begin{figure}[h]
\begin{center}
\scalebox{0.5}{\includegraphics{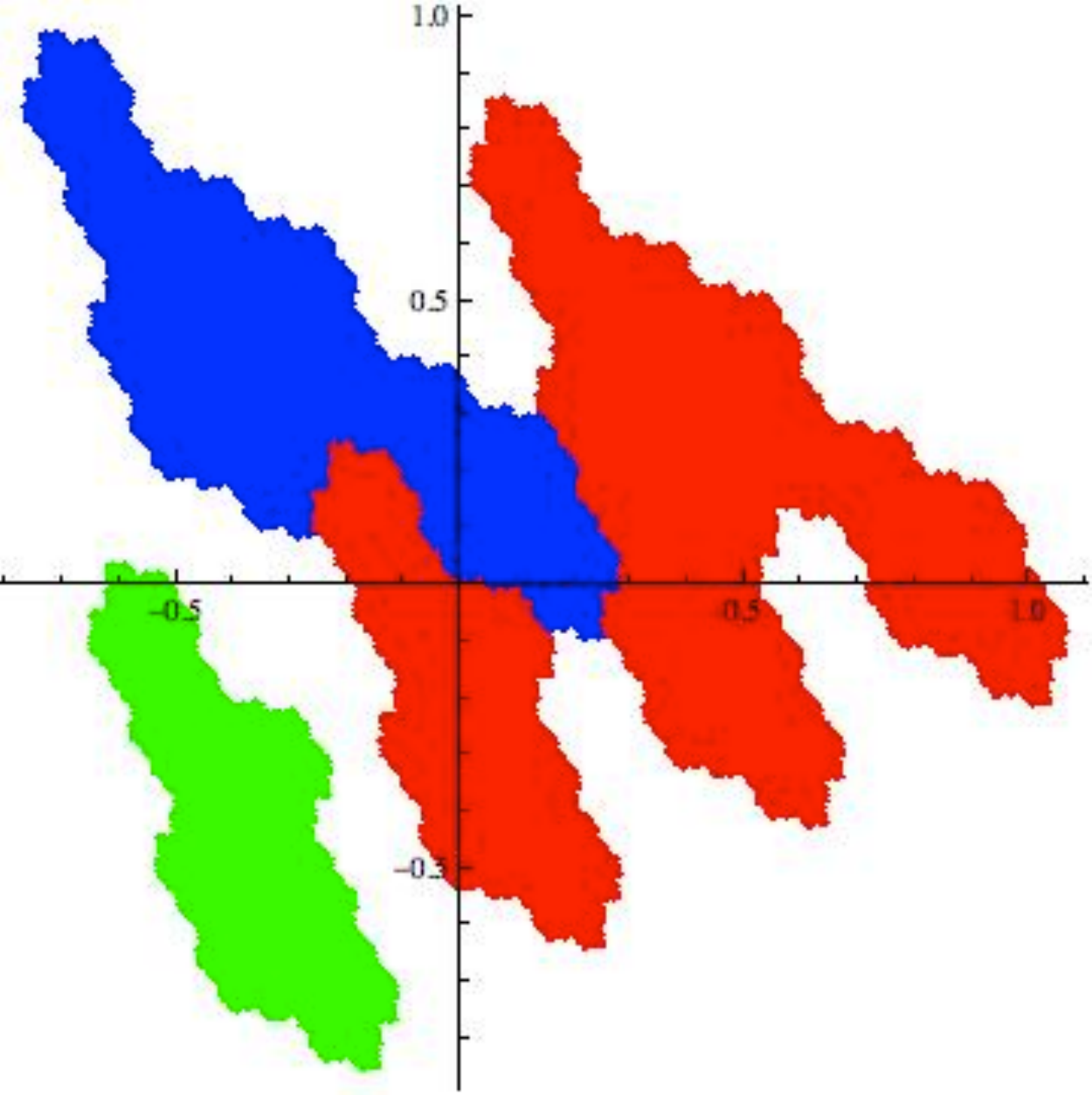}}
\caption{ The product Rauzy fractal of $\sigma_1$ and $\sigma_2$ with the sequence $u=222211111111111121$}
\end{center}
\end{figure}
\begin{figure}[h]
\begin{center}
\scalebox{0.5}{\includegraphics{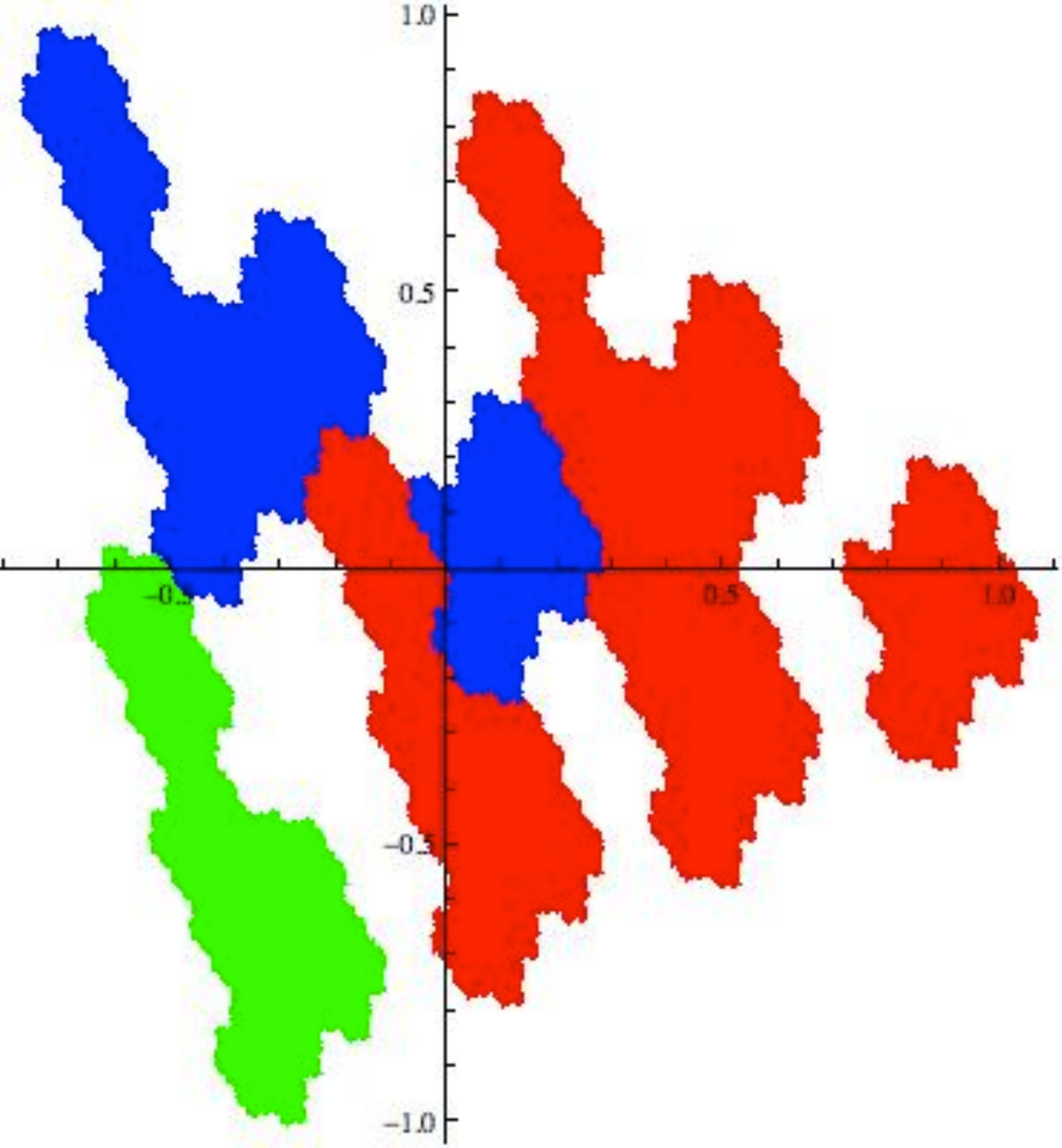}}
\caption{ The product  Rauzy  fractals of  $\sigma_1$ and $\sigma_2$ with the sequence $u =2222111111111111212$ }
\end{center}
\end{figure}
\begin{figure}[h]
\begin{center}
\scalebox{0.5}{\includegraphics{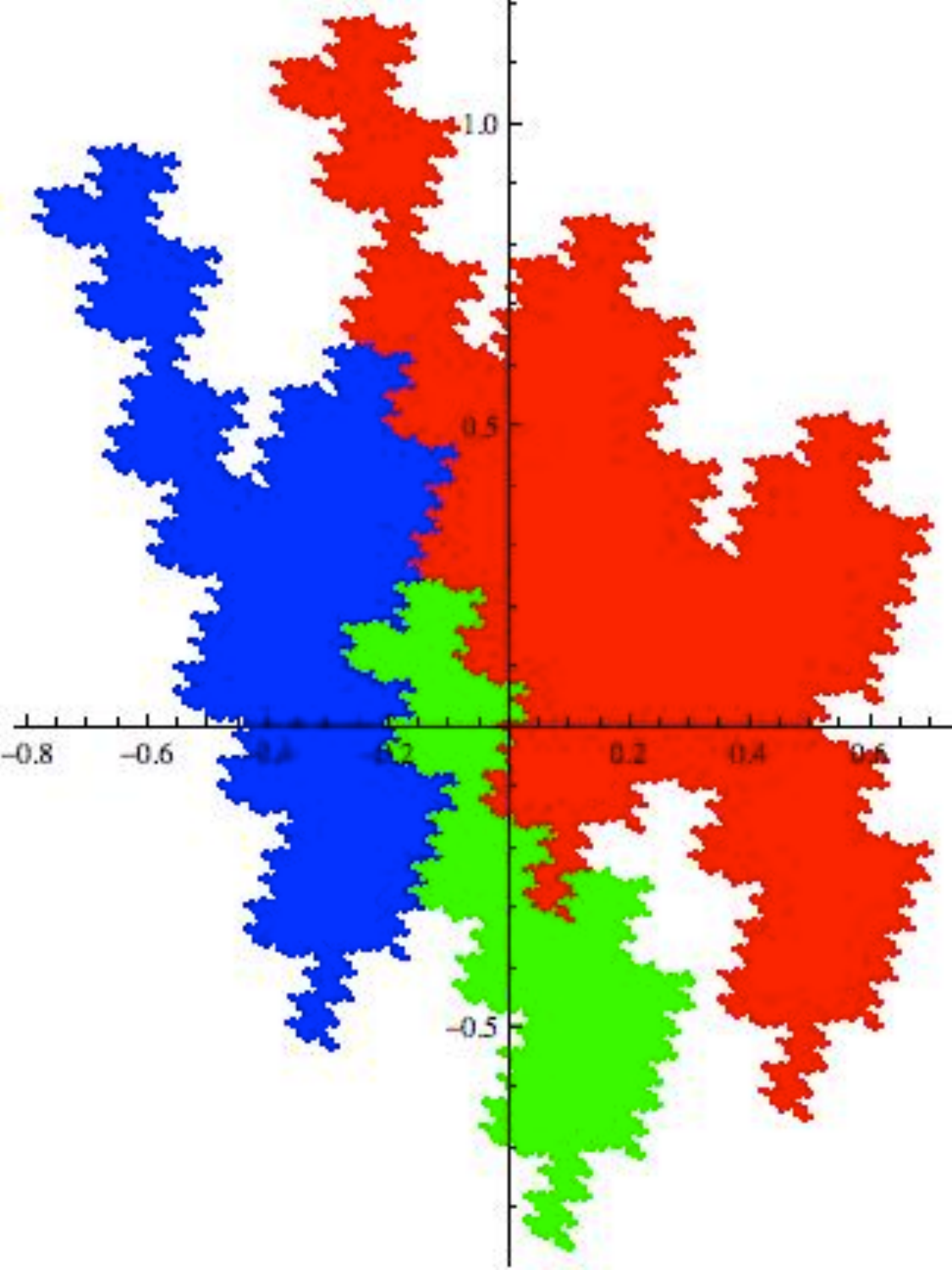}}
\caption{ The product Rauzy  fractals of  $\sigma_1$ and $\sigma_2$ with the sequence $u =1122112211221122$ }
\end{center}
\end{figure}
\section{additional remarks}
\subsection{ Substitutions with the same matrix}
For any $\sigma$, the set $\Rs$ contains a fundamental domain for the lattice $\Gamma$ acting on the contracting plane. More precisely, the union $\cup _{\gamma\in\Gamma}\Rs+\gamma$ is the contracting plane. The quotient $\Rs/\Gamma$ is a torus, and the map from the limit point to $\Gamma$ extends to a continuous map from the adic system $\Omega_{\sigma}$ to this torus, which is a semi-conjugacy between the shift on $\Omega_{\sigma}$ and a torus rotation.

On would of course like to prove that, if the matrix is unimodular, this semi-conjugacy is a measurable conjugacy (in the non-unimodular case, there must also be some $p$-adic component).

It seems very likely that the dimension of the boundary of $\R_{\sigma}$ is a measurable function of $\sigma$ which is shift invariant; hence, by ergodicity of the shift for any Bernoulli measure, it must be almost everywhere constant for this measure. However, this dimension is not constant, since it is generally not the same for the various substitutions; one would like to be able to compute the generic value for a given Bernoulli measure.

\subsection{The general adic case}

It would be much more interesting to extend this to the case of substitutions with different matrices, but this poses difficult problems.

We can define the symbolic system for any primitive sequence of substitutions, but after that we meet several obstructions. 

The first one is to find the equivalent of the Perron-Frobenius eigenvector, or the asymptotic direction for the stepped line. This is not always possible: the existence of this asymptotic direction is equivalent to the existence of a well-defined frequency of the letters, or to the unique ergodicity. But there are $S$-adic systems which are not uniquely ergodic, and we can even find such systems on 4 letters which are generated by a primitive sequence of unimodular substitutions.

However, such examples seem exceptional, and it might well be possible to find quite general conditions which ensure unique ergodicity of the $S$-adic system.

But if this preliminary step is satisfied, we still need to find the equivalent of the Pisot condition. One should not expect to find an exact equivalent: the Pisot condition only makes sense for algebraic numbers, and for most sequences of substitution, the equivalent of the Perron Frobenius eigenvector should have transcendant coordinates, and there is no equivalent for eigenvalues.

The best condition is probably to ask for a bounded distance between the stepped line associated with the limit point and the asymptotic line. Another equivalent condition is the condition of bounded balance; we say that a sequence is $C$-balanced if there exists a constant $C$ such that, for any two words $U,V$ of same length occuring in the sequence, and any letter $a$, $\left| |U|_a-|V|_a\right|\le C$.

Any sequence with bounded balance can define a generalized Rauzy fractal, and it would be very interesting to give sufficiently general condition which imply bounded balance; as shown by the case of unbalanced episturmian sequences, this does not seem to be an easy question.

\end{document}